\documentclass[11pt]{article}
\usepackage{graphicx}
\usepackage{amsmath,latexsym,amsbsy,amssymb, color}
\usepackage{graphicx}
\def\ds{\displaystyle}
\def\forall{\hbox{for all}~}
\def\L{{\bf L}}

\def\ve{\varepsilon}

\def\G{{\cal G}}

\def\R{I\!\!R}

\def\vp{\varphi}

\def\vs{\vskip 2em}

\def\v{\vskip 1em}

\def\O{{\cal O}}

\def\M{{\cal M}}
\def\S{{\cal S}}
\def\C{{\cal C}}

\def\ov{\overline}

\def\Tilde{\widetilde}

\def\bega{\begin{array}}
\def\enda{\end{array}}
\def\begi{\begin{itemize}}
\def\endi{\end{itemize}}

\def\bel{\begin{equation}\label}
\def\eeq{\end{equation}}
\def\sqr#1#2{\vbox{\hrule height .#2pt
\hbox{\vrule width .#2pt height #1pt \kern #1pt
\vrule width .#2pt}\hrule height .#2pt }}
\def\square{\sqr74}
\def\endproof{\hphantom{MM}\hfill\llap{$\square$}\goodbreak}

\newtheorem{theorem}{Theorem}[section]

\newtheorem{lemma}{Lemma}[section]
\newtheorem{proposition}{Proposition}[section]
\newtheorem{remark}{Remark}[section]
\newtheorem{definition}{Definition}[section]

\begin{document}

\title{\bf Generic Solutions to Controlled Balance Laws}
\vs
\author{Alberto Bressan$^{(*)}$ and 
Khai T.~Nguyen$^{(**)}$\\~~\\
 {\small $^{(*)}$~Department of Mathematics, Penn State University,}
 \\ {\small $^{(**)}$~Department of Mathematics, North Carolina State University.} \\~~\\
{\small E-mails: axb62@psu.edu,~khai@math.ncsu.edu.}
}
\maketitle

\begin{abstract} 
The paper is concerned with  a scalar balance law, where the source term depends on a control function $\alpha(t)$. 
Given a control $\alpha\in \L^\infty\bigl([0,T]\bigr)$, it is proved that, 
for  generic initial data $\bar u \in \C^3(\R)$, 
the solution has finitely many shocks, interacting at most two at a time. 
Moreover, at the terminal time $T$ no shock interaction occurs, and no new shock is formed.

In addition, a family of  optimal control problems is considered, including a running cost 
and a terminal cost.  An example is constructed where the optimal solution contains two
shocks merging exactly at the terminal time $T$. Such behavior persists under any suitably small perturbation of the flux, source, and cost functions, and of the initial data.
This shows that generic solutions of optimization problems have
different qualitative properties, compared with generic solutions to Cauchy problems.
\end{abstract}

\v
\section{Introduction}
\setcounter{equation}{0}
\label{s:1}

Consider the Cauchy problem for a scalar conservation law in one space dimension
with smooth, strictly convex flux
\bel{claw1} u_t+f(u)_x~=~0,\qquad\qquad  u(0,x)\,=\, \bar u(x),\qquad t\in [0,T].\eeq
For a generic initial data $\bar u\in \C^3(\R)$, it is well known that
the solution $u(t,\cdot)$ is piecewise continuous with a locally  finite number of shocks, 
which interact two at a time \cite{GS, Gk, S}.  
We recall that, on a complete metric space ${\cal Z}$, a property $P$ is {\bf generic} 
if it is satisfied by all  $z\in {\cal Z}'$, where ${\cal Z}'\subset {\cal Z}$
is a $\G_\delta$ set, i.e., it is the intersection of countably many open dense subsets in $\mathcal{Z}$.

Aim of  the present paper is to study regularity properties for generic solutions
to a controlled balance law, where the control $\alpha:[0,T]\mapsto \R^m$ acts on the source term
\cite{BM1, BS, U1, U2, U3}
\bel{cblaw}u_t+f(u)_x~=~g(t, x,u,\alpha),\qquad\qquad  u(0,x)\,=\, \bar u(x),\qquad t\in [0,T].\eeq
As basic assumptions, we consider
\begi
\item[{\bf (A1)} ]{\it The flux function $f$ is smooth and strictly convex, so that
$f''(u)\geq c_0>0$ for all $u\in\R$.  The source $g:[0,T]\times \R\times \R\times\R^m\mapsto\R$ is a bounded, smooth function, 
while $\alpha\in \L^\infty\bigl([0,T]\,;~\R^m\bigr)$ is a bounded, measurable control function.
}
\endi

Our first results show that the generic regularity properties proved in \cite{S, GS} remain valid
in the presence of a source term depending smoothly on $x$ but only measurably on $t$.
 
\begin{theorem}\label{t:1} Let $T>0$ and let $f,g,\alpha$ be given, satisfying {\bf (A1)}.  
Then there is a 
$\G_\delta$ set of initial data $\M\subset\C^3(\R)$ with the following property.
For every initial data $\bar u\in \M$, the 
solution to (\ref{cblaw}) is piecewise continuous, containing  finitely many shocks on any bounded domain $\bigl\{(t,x)\,;~t\in [0,T],~ x\in[-n,n]\bigr\}$. 
\end{theorem}

\begin{theorem}\label{t:2}  In the same setting of Theorem~\ref{t:1}, there exists a
$\G_\delta$ set $\M'\subseteq\M $ such that , for every initial data $\bar u\in \M'$, the 
solution to (\ref{cblaw}) has the following additional properties:
\begi
\item[(i)] The gradient $u_x$ blows up only at points where a new shock is formed.
\item[(ii)] Shocks interact only two at a time.
\item[(iii)] No shock interaction occurs at the terminal time~$T$.
\item[(iv)] No gradient blow-up occurs at the terminal time~$T$.

\endi
\end{theorem}

\begin{figure}[ht]
\centerline{\hbox{\includegraphics[width=15cm]{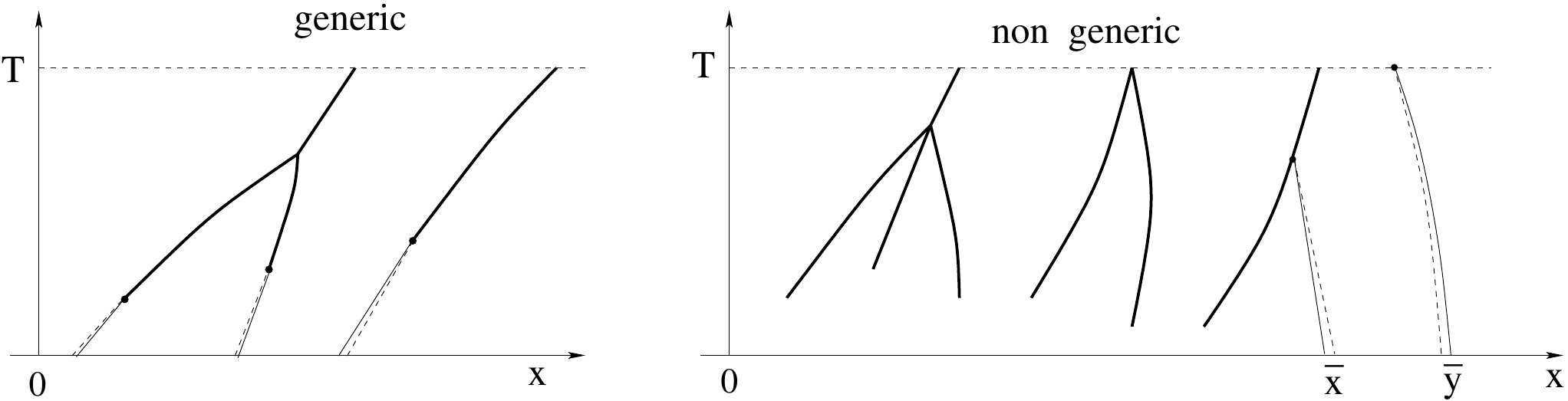}}}
\caption{\small  Left: a generic, structurally stable shock configuration.
Right: examples of non-generic configurations. These occur because:
(i) Three shocks interact simultaneously. (ii) Two shocks interact exactly at the terminal time $t=T$.  (iii) Along the characteristic starting at $\bar x$, the gradient blows up 
exactly at the time when another shock is reached. (iv) Along the characteristic starting at $\bar y$, the gradient blows up exactly at the terminal time $t=T$. 
These four features are all structurally 
unstable, since they no longer hold for an arbitrarily small perturbation of the initial data. }
\label{f:gen66}
\end{figure}

More precisely, the above statement (i) rules out the situation illustrated in Fig.~\ref{f:gen66}, right,
where the gradient $u_x\bigl(t, x(t)\bigr)\to -\infty$ as $t\to \tau-$, at the same time $\tau$
when the characteristic  $t\mapsto x(t)$ starting at $\bar x$ hits another shock, formed earlier.  In this case, 
no new shock is produced.
\v

Next, together with (\ref{cblaw}), we consider an optimal control problem of the form
\bel{min1}\hbox{minimize:} \qquad J[\alpha]~\doteq~\int_0^T\int \Phi\bigl(t,x,u(t,x), \alpha(t) \bigr)\, 
dx\, dt + \int \Psi\bigl(x,u(T,x)\bigr)\, 
dx\,.\eeq
Here $u=u(t,x)$ denotes the solution to (\ref{cblaw}) corresponding to the control $\alpha(\cdot)$, while $\Phi,\Psi$ are smooth functions, accounting for a running cost and a terminal cost, respectively.
In connection with this minimization problem, it is natural to speculate whether the 
optimal solution $u=u(t,x)$ shares the same regularity properties of generic solutions of 
(\ref{cblaw}), or else the minimization of some functional of the type (\ref{min1}) imposes distinctively particular properties.
The following questions thus arise:
\begi
\item[{\bf (Q1)}] {\it Is it true that, 
for a generic set of data $(f,g,\Phi,\Psi,\bar u)$,  the  optimal control problem (\ref{cblaw})-(\ref{min1}) 
has a
solution $u=u(t,x)$ which is piecewise continuous with finitely many shocks?}
\item[{\bf (Q2)}] {\it In a generic setting, is it true that, in a  optimal solution,  
shocks interact only two at a time and
no interaction occurs at the terminal time $T$?}
\endi

We leave {\bf (Q1)} as an open question.  
A negative answer to {\bf (Q2)} will be provided in Section~\ref{s:4}, namely
\begin{proposition}\label{p:11}
There exists a set of functions $(f,g,\Phi,\Psi,\bar u)$, open in the $\C^3$ topology, for which the solution to the optimization
problem (\ref{cblaw})-(\ref{min1})  contains precisely two shocks, interacting at the terminal  
time $t=T$.
\end{proposition}

We remark that, in connection with a minimization problem such as (\ref{cblaw})-(\ref{min1}), 
necessary conditions for optimality
 were obtained in \cite{BM1, BS}, in the form of a Pontryagin maximum principle.   
 These results  apply more generally to  $n\times n$ hyperbolic systems of balance laws.
However, to implement this approach one needs to assume
that the optimal solution is 
``shift differentiable" in the sense of \cite{BG, BL, BM2}. In particular, it should have isolated shocks interacting two at a time,
and no interaction should occur at the terminal time $T$.

The situation is somewhat different in the scalar case.   Indeed, 
by further developing the approach based on shift differentials, a remarkable result proved in 
\cite{U2}  shows that the mapping ``control-to-terminal cost"
$$\alpha(\cdot)~\mapsto~\int \Psi\bigl( x, u(T,x)\bigr)\, dx$$
can be differentiable under the only assumption that the terminal profile 
$u(T,\cdot)$ contains finitely many discontinuities, and no shock interaction occurs
 at the terminal time $T$.
In particular, no assumption about the structure of the BV solution 
(such as: only  pairwise shock interactions) is required at intermediate times 
$0<t<T$.

In view of Proposition~\ref{p:11}, however, there is a fairly large class of problems which 
do not fall within the scope of 
results  in \cite{BM1, BS, U2, U3}.  To treat these cases, additional necessary conditions 
for optimality would be needed.

A bit more promising is the case where only a running cost is present.
Indeed, if {\bf (Q1)} has a positive answer, 
and a generic optimal solution has finitely many shocks,
the analysis in \cite{U2} indicates that the mapping 
$$\alpha(\cdot)~\mapsto~\int_0^T \int \Phi\bigl(t,x, u(T,x), \alpha(t)\bigr)\,dx\, dt$$
would be differentiable also 
at the point  $\alpha^{opt}(\cdot)$ yielding the optimal control.

In the remainder of the paper, Section~\ref{s:2} contains a proof of Theorem~\ref{t:1}, 
while the proof of Theorem~\ref{t:2} is given in Section~\ref{s:3}.
Finally, a specific minimization problem, where the optimal solution has the behavior
described in Proposition~\ref{p:11}, is constructed in Section~\ref{s:4}.

\section{Proof of Theorem~\ref{t:1}}
\setcounter{equation}{0}
\label{s:2} 
The proof will be given in several steps.

{\bf 1.} Given  smooth functions $f,g$ and a measurable control function $\alpha $ as in {\bf (A1)}, for every initial
data $\bar{u}\in\mathcal{C}^3(\R) $ and $y\in\R$, we denote by 
$t\mapsto \bigl(\xi(t,y),v(t,y)\bigr)$ the solution to the system of ODEs for characteristics
\bel{ODE1}
\begin{cases}
\dot{\xi}~=~f'(v),\\[3mm]
\dot{v}~=~g(t,\xi,v,\alpha),
\end{cases}
\qquad\mathrm{with}\qquad\quad
\begin{cases}
\xi(0)~=~y,\\[3mm]
v(0)~=~\bar{u}(y).
\end{cases}
\eeq
Here and in the sequel, the upper dot denotes a derivative w.r.t.~time.
As long as the partial derivative $\xi_y(t,y)$ remains positive, the solution is thus implicitly 
determined by
\bel{sol}u\bigl(t, \xi(t,y)\bigr)~=~v(t,y).\eeq
Differentiating (\ref{ODE1}) w.r.t.~the initial point $y$, 
we obtain the system
\bel{ODE2}
\begin{cases}
\dot\xi_y ~=~f''(v) v_y\,,\\[3mm]
\dot v_y~=~g_\xi(t,\xi,v,\alpha)\xi_y + g_v(t,\xi,v,\alpha)\, v_y\,,
\end{cases}
\qquad\mathrm{with}\qquad\quad
\begin{cases}
\xi_y(0)~=~1,\\[3mm]
v_y(0)~=~\bar u_x(y).\end{cases}
\eeq
Call $T(y)$ the first time when the derivative $\xi_y(t,y)$ vanishes.
It is well known that new shocks are formed at points $\xi\bigl(T(y),y\bigr)$ where the function
$T(\cdot)$ attains a local minimum.
To prove that a solution contains finitely many shocks on bounded sets, we thus need to 
show that such points of local minima are isolated.
\v
{\bf 2.} 
The function $T(\cdot)$ is implicitly defined by 
\bel{Tid} \xi_y\bigl(T(y), y\bigr)~=~0.\eeq
We claim that, if (\ref{Tid}) holds, then
\bel{tey}\partial_t \xi_y\bigl(T(y), y\bigr)\,\not=\,0.\eeq
Indeed, since $f''(v)>0$, by the first equation in  (\ref{ODE2}) it suffices to show that $v_y\bigl(T(y), y\bigr)\not= 0$.
Assume, on the contrary, that 
\bel{zerod}v_y\bigl(T(y), y\bigr)~=~\xi_y\bigl(T(y), y\bigr)~=~0.\eeq
Since the functions $f'', g_\xi, g_v$ are locally bounded, by (\ref{ODE2}) it follows
$$ \bigl|\partial_t\,\xi_y(t)\bigr|~\leq~C \,\bigl|v_y(t,y)\bigr|,
\qquad\qquad \bigl|\partial_t \,v_y(t,y)\bigr|~\leq~C \,\bigl|\xi_y(t,y)\bigr|+C\bigl|v_y(t,y)\bigr|,$$
for some constant $C$.
In view of  (\ref{zerod}), by Gronwall's lemma we obtain
$$\xi_y(t,y)~=~v_y(t,y)~=~0\qquad\forall t\in \bigl[0,T(y)\bigr],$$
contradicting the assumption $\xi_y(0,y)=1$.  This proves our claim (\ref{tey}).

As a consequence,  if $\xi_y\bigl(T(y_0), y_0\bigr)=0$, by the implicit function theorem
the map $y\mapsto T(y)$ is well defined on a whole neighborhood of $y_0$. 
\v
{\bf 3.} By further differentiations, from (\ref{ODE2}) we obtain systems of ODEs determining the
higher derivatives $\partial^k_y \xi(t,y)$, $\partial_y^k v(t,y)$.  By the assumptions {\bf (A1)},
it follows that all maps 
$$y\mapsto \partial^k_y \xi(t,y),\qquad\qquad y\mapsto \partial_y^k v(t,y),$$
are smooth. On the other hand, since the control $\alpha(\cdot)$ is bounded and measurable, 
the maps
$t\mapsto \partial^k_y v(t,y)$ are Lipschitz continuous while the maps 
$t\mapsto \partial^k_y \xi(t,y)$ have $\C^{1,1}$ regularity (i.e., continuously differentiable with 
Lipschitz derivative).

 For notational convenience, for all $ (t,y)\in [0,\infty[\,\times\R$ we define 
\bel{theta}
\theta(t,y)\,\doteq\, \xi_{y}(t,y).
\eeq
At every  point $y_0\in \R$ where $T(\cdot)$ attains a local minimum, one has
\bel{Tt0}
\partial_{y}\theta\bigl(T(y_0),y_0\bigr)\,=\,0.
\eeq
If $\partial_{yy}\theta(T(y_0),y_0)\neq 0$, 
by  the implicit function theorem we can solve the equation
$$\partial_y\theta(t,y)\,=\,0,$$
and obtain a $\C^1$ function $y=\Tilde y(t)$, defined in a neighborhood of $t_0\doteq T(y_0)$,
such that
$$\theta_y\bigl( t, \Tilde y(t)\bigr)~=~0,\qquad\forall \quad t\in \,]t_0-\delta, \,t_0+\delta[\,,$$
for some $\delta>0$.   
Recalling (\ref{Tt0}), it is clear that the graphs of the two functions implicitly defined
by $\theta(t,y)=0$ and $\theta_y(t,y)=0$ both contain the point $\bigl( T(y_0), y_0)$, 
but have no other intersection in a neighborhood of such point.
This rules out the existence of any other local minimum for the
function $T(\cdot)$  in a neighborhood of $y_0$. 

We conclude that  all the local minima of the function $T(\cdot)$ are isolated provided that 
the system of three equations in two variables
\bel{cond1}\theta(t,y)\,=\,0,\qquad 
\partial_{y}\theta(t,y)\,=\,0,\qquad \partial_{yy}\theta(t,y)\,=\,0,
\eeq
has no solution.
\medskip

\begin{figure}[ht]
\centerline{\hbox{\includegraphics[width=8cm]{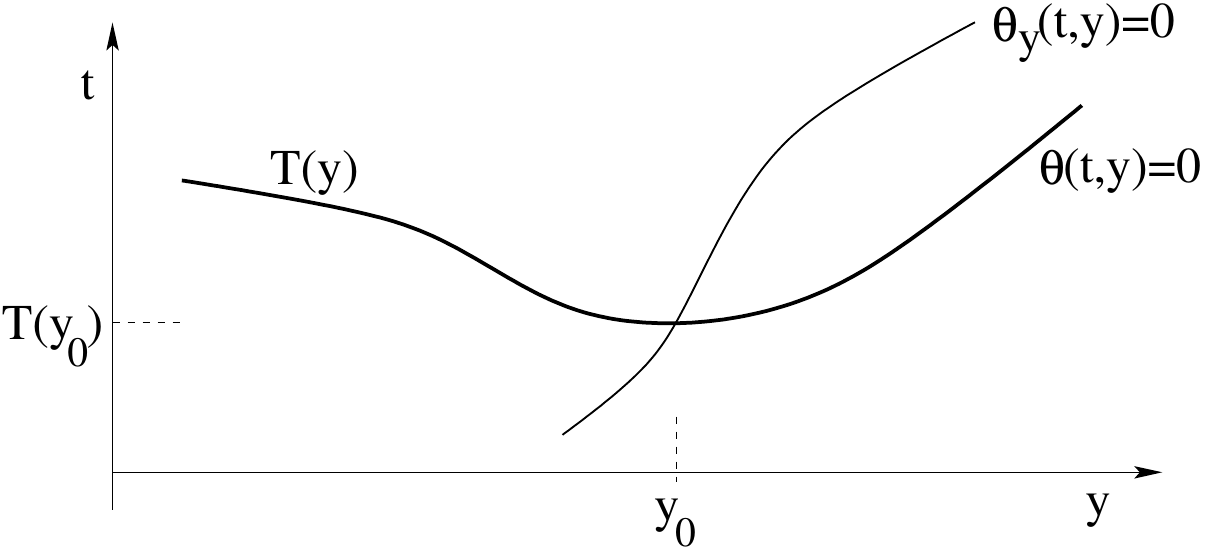}}}
\caption{\small  By the relative position of the curves implicitly defined by $\theta=0$ and $\theta_y=0$, one concludes that the local minimum of the function $T(\cdot)$ is isolated.}
\label{f:gen78}
\end{figure}

{\bf 4.} Observing that the functions $\theta, \theta_y, \theta_{yy}$ depend
continuously on the initial data $\bar u$ and on its first three derivatives, we conclude that for every $k\geq 1$ the set 
\bel{M-k}
\mathcal{M}_k~\doteq~\Big\{\bar{u}\in \mathcal{C}^3(\R)\,;~\hbox{the system (\ref{cond1}) has no solution for} ~(t,y)\in[0,T]\times [-k,k]\Big\}
\eeq
 is open in $\mathcal{C}^3(\R)$. 
 
 In the next steps, using standard tools
from differential geometry~
\cite{Bloom, GG, GS, Gk},
we will show that  each $\mathcal{M}_k$  is dense in $\C^3(\R)$.

%
%
\v

{\bf 5.} Consider any $\bar u\in \C^3(\R)$. 
For any given $\bar{y}\in \R$,  we define the 3-parameter 
family of perturbed initial data $\bar{u}^\lambda \in\mathcal{C}^3(\R)$,
$\lambda=(\lambda_1,\lambda_2,\lambda_3)\in [-1,1]^3$, by setting 
\bel{ula}\bar {u}^\lambda (y)~\doteq~\bar{u}(y) +\eta(y-\bar{y})\cdot \left(\lambda_1\cdot (y-\bar y)+\lambda_2\cdot {(y-\bar{y})^2\over 2}+\lambda_3\cdot {(y-\bar{y})^3\over 6}\right).\eeq
Here $\eta:\R\mapsto [0,1]$ is a smooth cutoff function such that 
\bel{cutoff}\eta(y)~=~\left\{ \bega{rl} 1\quad &\hbox{if}\quad |y|\leq 1,\\[1mm]
0\quad &\hbox{if}\quad |y|\geq 2.\enda\right.\eeq
The partial derivatives of $\bar u^\lambda$ at the point $\bar y$ 
are computed by
\[
\begin{cases}
\partial_{\lambda_1}\bar u^\lambda(\bar{y})~=~0,
\\[2mm]
\partial_{\lambda_2}\bar{u}^\lambda(\bar{y})~=~\partial_{\lambda_2y}
\bar{u}^\lambda(\bar{y})~=~0,\\[2mm]
\partial_{\lambda_3}\bar{u}^\lambda(\bar{y})~=~\partial_{\lambda_3y}
\bar{u}^\lambda(\bar{y})~=~\partial_{\lambda_3yy}\bar{u}^\lambda(\bar{y})~=~0,
\end{cases}
\qquad\quad
\begin{cases} \partial_{\lambda_1y}
\bar u^\lambda(\bar{y})~=~1,\\[2mm]\partial_{\lambda_2yy}\bar{u}^\lambda(\bar{y})~=~1,\\[2mm]
\partial_{\lambda_3yyy}\bar{u}^\lambda(\bar{y})~=~1.
\end{cases}
\]
In the following,  we denote by $\xi^\lambda, v^\lambda$ the corresponding solutions to (\ref{ODE1}), with $\bar u$ replaced by $\bar u^\lambda$.  Recalling (\ref{theta}),  we also 
write $\theta^\lambda(t,y) \doteq \xi^\lambda_y(t,y)$.
At time $t=0$ we have 
\bel{inc}
\left\{\bega{ll}
\partial_{\lambda_i}\xi^{\lambda}(0,\bar{y})~=~\partial_{\lambda_i}\theta^{\lambda}(0,\bar{y})~=~0,
\qquad i=1,2,3,&\\[2mm]
\partial_{\lambda_1}v^{\lambda}(0,\bar{y})~=~0,\quad &\partial_{\lambda_1}v^{\lambda}_y(0,\bar{y})~=~1,\\[2mm]
\partial_{\lambda_2}v^{\lambda}(0,\bar{y})~=~\partial_{\lambda_2}v^{\lambda}_{y}(0,\bar{y})
~=~0,\quad &\partial_{\lambda_2}v^{\lambda}_{yy}(0,\bar{y})~=~1,\\[2mm]
\partial_{\lambda_3}v^{\lambda}(0,\bar{y})~=~\partial_{\lambda_3}v^{\lambda}_{y}(0,\bar{y})~=~\partial_{\lambda_3}v^{\lambda}_{yy}(0,\bar{y})~=~0,\quad &\partial_{\lambda_3}v^{\lambda}_{yyy}(0,\bar{y})~=~1.\enda
\right.
\eeq
The Cauchy problem  for $(\xi^{\lambda}, v^{\lambda},\theta^{\lambda})$ takes the form
\bel{E-theta}
\begin{cases}
\ds{\partial \over\partial t}\xi^{\lambda}(t,y)~=~f'\bigl(v^{\lambda}(t,y)\bigr),\\[4mm]
\ds{\partial \over \partial t}v^{\lambda}(t,y)~=~g\bigl(t,\xi^{\lambda},v^{\lambda},\alpha(t)\bigr),\\[4mm]
\ds{\partial \over \partial t}\theta^{\lambda}(t,y)~=~f''\bigl(v^{\lambda}(t,y)\bigr)\cdot v^{\lambda}_y(t,y),
\end{cases}
\qquad
\begin{cases}
\xi^{\lambda}(0,y)~=~y,\\[4mm]
v^{\lambda}(0,y)~=~\bar {u}^\lambda(y),\\[4mm]
\theta^{\lambda}(0,y)~=~1.
\end{cases}
\eeq
For notational convenience, we introduce the functions
\[\left\{ \bega{rl} 
a(t)&\doteq~f''\bigl(v^{\lambda}(t,\bar{y})\bigr)~=~{f''\bigl(v(t,\bar{y})\bigr)},\\[1mm]
b(t)&\doteq~g_{\xi}\bigl(t,\xi^{\lambda}(t,\bar{y}),v^{\lambda}(t,\bar{y}),\alpha(t)\bigr)~=~{g_{\xi}\bigl(t,\xi(t,\bar{y}),v(t,\bar{y}),\alpha(t)\bigr)},\\[1mm]
c(t)&\doteq~g_{v}\bigl(t,\xi^{\lambda}(t,\bar{y}),v^{\lambda}(t,\bar{y}),\alpha(t)\bigr)~=~{g_{v}\bigl(t,\xi(t,\bar{y}),v(t,\bar{y}),\alpha(t)\bigr)}.
\enda\right.
\]
Differentiating the first two equations in (\ref{E-theta}) w.r.t $\lambda_i$, $i=1,2,3$, one obtains
\bel{E2}
\begin{cases}
\ds{\partial \over\partial t}\partial_{\lambda_i}\xi^{\lambda}(t,y)~=~f''(v^{\lambda})\cdot\partial_{\lambda_i}v^{\lambda}(t,y),\\[4mm]
\ds{\partial \over\partial t}\partial_{\lambda_i}v^{\lambda}(t,y)~=~g_{\xi}\cdot \partial_{\lambda_i}\xi^{\lambda}(t,y)+g_{v}\cdot \partial_{\lambda_i}v^{\lambda}(t,y),\\[4mm]
\end{cases}
\quad
\begin{cases}
\partial_{\lambda_i}\xi^{\lambda}(0,\bar{y})~=~0,\\[4mm]
\partial_{\lambda_i}v^{\lambda}(0,\bar{y})~=~0.
\end{cases}
\eeq
This yields
\bel{sc}
\partial_{\lambda_i}\xi^{\lambda}(t,\bar{y})~=~\partial_{\lambda_i}v^{\lambda}(t,\bar{y})~=~0\qquad\forall t\in [0,T].
\eeq
Differentiating (\ref{E2})  w.r.t.~$y$ and observing that some terms vanish because of (\ref{sc}),
we check that the maps  $t\mapsto \bigl(\partial_{\lambda_i}\theta^{\lambda}(t,\bar{y}),~\partial_{\lambda_i}v^{\lambda}_y(t,\bar{y})\bigr)$ provide
a solution to  the system of ODEs
\bel{C-ODE}
\begin{cases}
\dot X(t)~=~a(t)\cdot Y(t),\\[2mm]
\dot Y(t)~=~b(t)\cdot X(t)+c(t)\cdot Y(t),
\end{cases}
\eeq
%
with initial data
\[ 
X(0)\,=\,\partial_{\lambda_1}\theta^{\lambda}(0,\bar{y})\,=\,0,\qquad 
Y(0)\,=\,\partial_{\lambda_1}v^{\lambda}_{y}(0,\bar{y})\,=\,1,\qquad \hbox{for}~i=1,
\]
\[ X(0)\,=\,
\partial_{\lambda_i}\theta^{\lambda}(0,\bar{y})\,=\,0,\qquad Y(0)\,=\,\partial_{\lambda_i}v^{\lambda}_{y}(0,\bar{y})\,=\,0,\qquad\hbox{for}~ i=2,3.
\]
This implies
\bel{sc1}
\partial_{\lambda_i}\theta^{\lambda}(t,\bar{y})~=~\partial_{\lambda_i}v^{\lambda}_y(t,\bar{y})~=~0\qquad\forall t\in [0,T],~~ i=2,3.
\eeq
Differentiating  (\ref{E-theta}) once again w.r.t.~$y$ and using (\ref{sc1}) we obtain that, for $i=2,3$, the maps  
$t\mapsto \bigl(\partial_{\lambda_i}\theta^{\lambda}_y(t,\bar{y}),~\partial_{\lambda_i}v^{\lambda}_{yy}(t,\bar{y})\bigr)$ provide
a solution to  the same system of ODEs
(\ref{C-ODE}), with initial data
\[
X(0)\,=\,
\partial_{\lambda_2}\theta^{\lambda}_y(0,\bar{y})\,=\,0,\qquad 
Y(0)\,=\,\partial_{\lambda_2}v^{\lambda}_{yy}(0,\bar{y})\,=\,1,
\]
\[X(0)\,=\,
\partial_{\lambda_3}\theta^{\lambda}_y(0,\bar{y})\,=\,0,
\qquad Y(0)\,=\,\partial_{\lambda_3}v^{\lambda}_{yy}(0,\bar{y})\,=\,0,
\]
in the cases $i=2$ and $i=3$, respectively.
In particular, we have
\bel{sc2}
\partial_{\lambda_3}\theta^{\lambda}_y(t,\bar{y})~=~\partial_{\lambda_3}v^{\lambda}_{yy}(t,\bar{y})~=~0\qquad\forall t\in [0,T].
\eeq
Finally, differentiating (\ref{E-theta}) a third time w.r.t.~$y$ and using (\ref{sc2}), we check that the map $t\mapsto \bigl(\partial_{\lambda_3}\theta^{\lambda}_{yy}(t,\bar{y})\,,~\partial_{\lambda_3}v^{\lambda}_{yyy}(t,\bar{y})\bigr)$ provides a solution to (\ref{C-ODE}) with initial data
\bel{I-la-3}X(0)\,=\,
\partial_{\lambda_3}\theta^{\lambda}_{yy}(0,\bar{y})\,=\,0,\qquad 
Y(0)\,=\,\partial_{\lambda_3}v^{\lambda}_{yyy}(0,\bar{y})\,=\,1.
\eeq

{\bf 6.}
By the analysis in the previous step, calling $t\mapsto \bigl( \ov X(t), \ov Y(t)\bigr)$
the solution to  (\ref{C-ODE}) with initial data
\bel{ODE-id}\bigl( \ov X(0),\,\ov Y(0)\bigr)~=~(0,1),\eeq 
one has
\[
\partial_{\lambda_1}\theta^{\lambda}(t,\bar{y})~=~\partial_{\lambda_2}\theta^{\lambda}_{y}(t,\bar{y})~=~\partial_{\lambda_3}\theta^{\lambda}_{yy}(t,\bar{y})~=~\ov X(t),
\]
\[
\partial_{\lambda_1}v^{\lambda}_y(t,\bar{y})~=~\partial_{\lambda_2}v^{\lambda}_{yy}(t,\bar{y})~=~\partial_{\lambda_3}v^{\lambda}_{yyy}(t,\bar{y})~=~\ov Y(t).
\]
Recalling (\ref{sc1}) and (\ref{sc2}), we thus have
\bel{det}
\det\Big[D_{\lambda}\bigl(\theta^{\lambda}(t,\bar{y}),\,\partial_{y}\theta^{\lambda}(t,\bar{y}),\,\partial_{yy}\theta^{\lambda}(t,\bar{y})\bigr)\Big] ~=~\det \begin{bmatrix}
\ov X(t)&*&*\\[2mm]
0&\ov X(t)&*\\[2mm]
0&0&\ov X(t)
\end{bmatrix}
~=~\ov X^3(t).
\eeq
We now observe that , if $\ov X(\tau)=0$  at some time $\tau\geq 0$, then $\ov Y(\tau)\neq 0$. Otherwise
the solution to (\ref{C-ODE}) with data $\bigl( \ov X(\tau), \ov Y(\tau)\bigr)=(0,0)$
would be identically zero, against the assumption (\ref{ODE-id}).
%

On the other hand, by (\ref{ODE2}) the map 
$$t\,\mapsto \,\bigl(X(t), Y(t)\bigr)\,=\,\bigl(\theta^{\lambda}(t,\bar{y}), v^{\lambda}_y(t,\bar{y})\bigr)$$
provides another solution to 
(\ref{C-ODE}), with  
 $$\bigl(X(0), Y(0)\bigr)~=~\big(\theta^{\lambda}(0,\bar{y}), v^{\lambda}_y(0,\bar{y})\big)~=~(1,\, \bar u_x(\bar y)+\lambda_1).$$
 We claim that, for this second solution,  
\bel{g-c} X(\tau)~=~
\theta^{\lambda}(\tau,\bar{y})~\neq~0.
\eeq
Indeed,  if $X(\tau)=0$, since we are also assuming that $\ov X(\tau)=0$,  
by linearity it follows 
\[ X(t)~=~ {Y(\tau)\over \ov Y(\tau)} \, \ov X(t)
\qquad\forall t\geq 0.
\]
This  yields a contradiction, because it implies
\[
1~=~X(0)~=~\theta^{\lambda}(0,\bar{y})~\not=~{Y(\tau)\over \ov Y(\tau)}\,
 \ov X(0)~=~0.
\]
Combining (\ref{det}) with (\ref{g-c}), we thus conclude
\[
\big|\theta^{\lambda}(t,\bar{y})\big|+\bigg|\mathrm{det}\Big[D_{\lambda}\bigl(\theta^{\lambda}(t,\bar{y})\,,~\partial_{y}\theta^{\lambda}(t,\bar{y})\,,~\partial_{yy}\theta^{\lambda}(t,\bar{y})\bigr)\Big]\bigg|~>~0,\qquad\forall t\geq 0.
\]
In particular, by continuity there exists $r_{\bar{y}}>0$ such that 
\bel{f-R}
\big|\theta^{\lambda}(t,y)\big|+\bigg|\mathrm{det}D_{\lambda}\bigl(\theta^{\lambda}(t,y),\,\partial_{y}\theta^{\lambda}(t,y),\,\partial_{yy}\theta^{\lambda}(t,y)\bigr)\bigg|~>~0,
\eeq
for all $(t,y)\in R_{\bar{y}}=[0,T]\times [\bar{y}-r_{\bar{y}},\,\bar{y}+r_{\bar{y}}]$ and $\lambda\in [-1,1]^3$.
\v
{\bf 7.} Next, covering the compact set $[0,T]\times [-k,k]$ 
with finitely many rectangles $R_{y^{\ell}}=[0,T]\times\big[y^{\ell}-r_{y^{\ell}},y^{\ell}+r_{y^{\ell}}\big]$, $\ell=1,\cdots, N$, we consider a family of combined perturbations defined as follows. For $\lambda= (\lambda^{1},\cdots,\lambda^{N})\in [-1,1]^{3N}$, 
we set
\bel{ini-p}
\bar u^\lambda(y)~=~\bar{u}(y) +\sum_{\ell=1}^{N}\eta{(y-y^{\ell})}\cdot\left(\lambda^{\ell}_1\,(y-y^{\ell})+\lambda^{\ell}_2\, {(y-y^{\ell})^2\over 2}+\lambda^{\ell}_3\, {(y-y^{\ell})^3\over 6}\right).
\eeq
We then define  the $\mathcal{C}^3$ map $\Phi:[-1,1]^{3N}\times [0,T]\times\R\mapsto  \R^3$ by setting
\bel{phi}
\Phi(\lambda,t,y)~\doteq~\bigl(\theta^{\lambda}(t,y),\,\partial_{y}\theta^{\lambda}(t,y),\,\partial_{yy}\theta^{\lambda}(t,y)\bigr),\qquad (\lambda,t,y)\in [-1,1]^{3N}\times[0,T]\times\R.
\eeq
Notice that when $\lambda^{j}=0$ for some $j\in \{1,\ldots,N\}$,  the perturbation at $y^j$ is absent  from $\bar u^\lambda(y)$. Hence it does not  affect the solution $\theta^{\lambda}(t,y)$. Consequently, differentiating only w.r.t.~the components 
$(\lambda_1^\ell,\lambda_2^\ell, \lambda_3^\ell)$ for one particular $\ell$, at 
$\lambda={\bf 0}\in \R^{3N}$  in view of  (\ref{f-R})  we obtain
\[
\Big|\theta^{\lambda}(t,y)_{|_{\lambda={\bf 0}}}\Big|+{\bigg|\mathrm{det}\Big[D_{\lambda^{\ell}}(\theta^{\lambda}(t,y),\partial_{y}\theta^\lambda(t,y),\partial_{yy}\theta^\lambda(t,y))\Big]_{\lambda={\bf 0}}\bigg|}~>~0,\qquad\forall (t,y)\in R_{y^{\ell}}\,.
\]
Again by continuity, choosing $\ve>0$ sufficiently small, for all $(\lambda,t,y)\in [-\ve,\ve]^{3N}\times [0,T]\times [-k,k]$ we have 
\[
\hbox{either}\quad 
\Phi(\lambda,t,y)\neq 0\qquad\mathrm{or}\qquad \mathrm{rank} \bigl[D_{\lambda^\ell}\Phi(\lambda,t,y)\bigr]~=~3.
\]
By the transversality theorem \cite{Bloom, GG}, for a dense set of values $\lambda\in [-\ve,\ve]^{3N}$, the map $(t,y)\mapsto \Phi(\lambda, t,y)$ is transversal to the zero-dimensional submanifold
$\{{\bf 0}\}\subset\R^3$, restricted to the domain where $(t,y)\in [0,T]\times [-k,k]$. 
Of course, this can happen only if ${\bf 0} = (0,0,0)\in \R^3$ is not in the range
of $\Phi(\lambda, \cdot,\cdot)$.
We conclude that the set $\mathcal{M}_{k}$ defined at (\ref{M-k}) is dense in $\mathcal{C}^3(\R)$.

\v
{\bf 8.} Since every set $\M_k$ is open and dense, 
 the set 
\bel{Mint}
\mathcal{M}~\doteq~\bigcap_{k\geq 1}\mathcal{M}_k\eeq is a $\mathcal{G}_{\delta}$ subset of $\mathcal{C}^3(\R)$. 
We claim that,
for every $\bar{u}\in \mathcal{M}$, the  solution to (\ref{cblaw}) is piecewise continuous with finitely many shocks on any bounded domain 
$\bigl\{(t,x);\, t\in [0,T], ~x\in [-n,n] \bigr\}$.

Indeed, consider an initial data $\bar u\in \M$.   Since by assumption $g$ is bounded,
the corresponding solution $u=u(t,x)$  is uniformly bounded on $[0,T]\times \R$.
As a consequence, all wave  speeds remain uniformly bounded, say
$\bigl| f(u(t,x))\bigr|\leq \lambda^*$ for some upper bound $\lambda^*>0$.
Given $n\geq 1$, we now choose an integer $k$ large enough so that 
$$[ -n,n]~\subset~[-k + \lambda^*T\,,~k-\lambda^* T].$$
Since $\bar u\in \M_k$, this implies that the solution $u$ is piecewise continuous with finitely many shocks on $[0,T]\times [-n,n]$, proving the theorem.
\endproof

\begin{remark}\label{r:21}
{\rm Consider a point $P=\bigl(\tau, \xi(\tau,y)\bigr)$, 
along the characteristic starting at $y$,
where a new shock is formed.  Under the generic conditions
\bel{sf0}\xi_y(\tau,y)=0,\qquad\quad\xi_{yy}(\tau,y)=0,\quad\qquad \xi_{yyy}(\tau,y)\not= 0,\eeq
for $t>\tau$ the asymptotic behavior of the new shock near $P$ is the same as for 
a scalar conservation law, see  \cite{Gk}.    Notice that in (\ref{sf0}) one must have
$\xi_{yyy}>0$, because we are considering a local minimum of $\xi_y$.

We observe that, after a shock is formed, 
it can merge with other shocks but it cannot disappear in finite time.
Indeed, call $u^-(t), u^+(t)$ the left and right states across the shock,
and $u_x^-(t), u_x^+(t)$ the gradients of $u(t,\cdot)$ to the left and to the right of the shock. 
Denoting by $x(t)$ the position of the shock at time $t$, the  Rankine-Hugoniot 
condition yields
$$\dot x(t)~=~{f(u^-(t))- f(u^+(t))\over u^-(t)-u^+(t)}\,.$$
Since the flux $f$ is convex, the Lax admissibility conditions imply $u^+(t)< u^-(t)$.
For the balance law  (\ref{cblaw}) we  thus have
$$\bega{rl} \ds{d\over dt} \bigl(u^-(t)- u^+(t)\bigr)&=~g\bigl(t, x(t), u^-(t), \alpha(t)\bigr) -g\bigl(t, x(t), u^+(t), \alpha(t)\bigr) \\[3mm]
&\qquad\qquad - \Big( f'\bigl(u^-(t)\bigr) - \dot x(t)\Big) u_x^-(t) + \Big( f'\bigl(u^+(t)\bigr) - \dot x(t)\Big) u_x^+(t)\\[3mm]
&=~\O(1) \cdot\bigl(u^-(t)- u^+(t)\bigr).\enda$$
Here and throughout the sequel, the Landau symbol $\O(1)$ denotes a uniformly bounded quantity.  We conclude  that the size of the shock can decrease at exponential rate, but will not
become zero in finite time.

}
\end{remark}

\begin{remark}
{\rm If $\bar u\in\M$ is a generic initial data satisfying the conditions in Theorem~\ref{t:1},
there may exist two  characteristics starting at $y_1<y_2$ and a time $\tau>0$
such that
$$
\theta(\tau, y_1) \, =\,\theta_y(\tau, y_2)\,=\,0,\qquad 
\theta_y(\tau, y_1) \, =\,\theta_y(\tau, y_2)\,=\,0,
$$
with 
\bel{ngen1}
\xi(\tau, y_1) \, =\,\xi(\tau, y_2)~\doteq~x^*.\eeq
We observe that in this case there must be a shock, formed earlier, passing through
the point $(\tau, x^*)$.  Otherwise, all characteristics starting at points $y\in [y_1, y_2]$
would meet at the same point $(\tau, x^*)$.
This would imply $\xi_y(\tau,y) = \xi_{yy}(\tau,y)=\xi_{yyy}(\tau,y)=0$ 
for all $y\in [y_1, y_2]$,
contradicting the assumption that the solution with initial data $\bar u\in\M$ has generic behavior.
}
\end{remark}

\begin{figure}[ht]
\centerline{\hbox{\includegraphics[width=15cm]{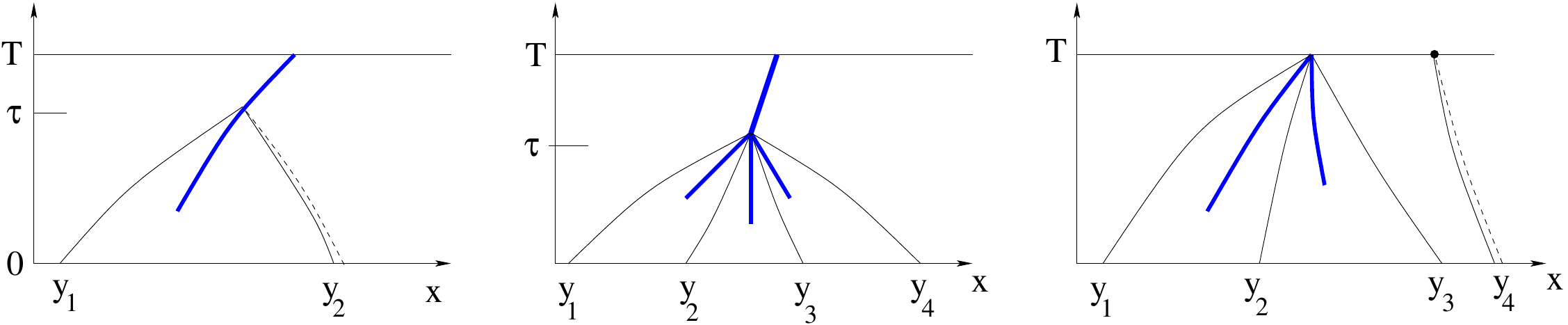}}}
\caption{\small Four non-generic configurations. In each case, the solution satisfies
a system where the number of equations is strictly larger than the number of variables. } 
\label{f:gen82}
\end{figure}

\section{Proof of Theorem~\ref{t:2}}
\setcounter{equation}{0}
\label{s:3} 

 Aim of this section we show that, for generic initial data $\bar u\in \M$, 
the four additional properties (i)--(iv) stated in Theorem~\ref{t:2} hold. 
The underlying motivation is that, if any one of these properties fails, then the 
solution $(\xi,v)$ of (\ref{ODE1}) satisfies a system  where the number of equations is strictly larger than the number of variables. More precisely:
\begi
\item If (i) fails (see Fig.~\ref{f:gen82}, left), then there exist two points
$y_1\not= y_2$ and a time $\tau\in \,[0,T]$  such that the following four identities holds:
\bel{1-1} 
\xi(\tau,y_1)\,=\,\xi(\tau, y_2),\qquad\quad \xi_y(\tau,y_2)\,=\,0,\qquad \quad
\xi_{yy}(\tau,y_2)\,=\,0,\eeq
\bel{1-2}
\bega{l} \ds\int_{y_1}^{y_{2}}\bar{u}(x)dx+\int_{0}^\tau\int^{\xi(t,y_2)}_{\xi(t,y_1)}g(t,x,u,\alpha)\,dxdt
+\int_{0}^\tau\bigl[f(v(t,y_1))-f'(v(t,y_1))v(t,y_1)\bigr]dt\\[4mm]
\qquad\qquad\ds -\int_{0}^\tau\bigl[f(v(t,y_2))-f'(v(t,y_2))v(t,y_2)\bigr]dt
~=~0.\enda
\eeq
W.l.o.g., in (\ref{1-2}) we assumed $y_1<y_2$. Notice that 
(\ref{1-2}) is obtained by integrating the balance law over the triangular region 
bounded by the two characteristics starting at $y_1$, $y_2$ and hitting the shock at time $\tau$.
This yields a system of  four equations for the three variables $y_1, y_2,\tau$.

\item If (ii) fails (see Fig.~\ref{f:gen82}, center), then there exist four points
$y_1<y_2<y_3<y_4$ and a time $\tau\in [0,T]$ such that, for $ i=1,2,3$, the following identities are satisfied:
\bel{2-1}
\xi(\tau, y_{i+1}) \,=\, \xi(\tau,y_i),\eeq
\bel{2-2}
\bega{l} \ds\int_{y_i}^{y_{i+1} }\bar{u}(x)dx+\int_{0}^\tau\int^{\xi(  t,y_{i+1})}_{\xi(  t,y_i)}g(t,x,u,\alpha)\,dxd  t
+\int_{0}^\tau\bigl[f(v(  t,y_i))-f'(v(  t,y_i))v(  t,y_i)\bigr]d  t\\[4mm]
\qquad\qquad\ds -\int_{0}^\tau\bigl[f(v(  t,y_{i+1}))-f'(v(  t,y_{i+1}))v(  t,y_{i+1})\bigr]d  t
~=~0.\enda
\eeq
This yields a system of  six equations for the five variables $y_1, y_2, y_3, y_4,\tau$.

\item If (iii) fails (see Fig.~\ref{f:gen82}, right), then there exist three points
$y_1<y_2<y_3$ such that, for $ i=1,2$, the identities (\ref{2-1}) and (\ref{2-2})
are satisfied with $\tau=T$. This yields a system of 
four equations for the three variables $y_1, y_2, y_3$.

\item If (iv) fails (see Fig.~\ref{f:gen82}, right), then there exists a point
$y_4$ such that the identities
\bel{y4d}\xi_y(T,y_4)\,=\,0,\qquad \quad
\xi_{yy}(T,y_4)\,=\,0,\eeq
are satisfied. This yields a system of 
two equations for the single variable $y_4$.
\endi

We remark, however, that a proof of Theorem~\ref{t:2} based on the transversality
theorem is not possible, because the double integrals in (\ref{1-2}) and (\ref{2-2}) 
include the values of the discontinuous solution $u=u(t,x)$, 
and may not depend smoothly on the initial data.
We will thus use a different argument, where the non-generic 
singularities are removed one at a time.

\v
Let $\M\subset \C^3(\R)$ be the $\G_\delta$-set defined at (\ref{Mint}). 
%
For every initial data $\bar u\in \M$, the solution $u=u(t,x)$ of (\ref{cblaw})
thus contains a locally finite number of shocks. 
For a fixed integer $\nu\geq 1$, let
\bel{SP} \S(\bar u)~\doteq~\Big\{ P_i=(t_i, x_i)\in [0,T]\times [-\nu,\nu],\quad i=1,\ldots,q\Big\},\eeq
 be the list of all points $(\bar t,\bar x)$ such that
\begi
\item[(i)] either  the relations (\ref{sf0}) 
holds at $\bigl(t,\,\xi(t,y)\bigr) =(\bar t,\bar x)$, so that at this point the gradient $u_x$ blows up,
\item[(ii)] or else  at  $(\bar t,\bar x)$ two or more shocks merge together.
\endi
According to Theorem~\ref{t:1}, if $\bar u\in \M$ then this set of points is finite. 
To each point $P_i\in \S(\bar u)$ we attach an integer $N_i$, describing by
how much the solution $u$ fails to be ``generic" at this point. 

\begin{definition}\label{def:31} Let $P= (\bar t,\bar x)\in\S $.   Assume that there
exist $m$ distinct
characteristics, starting at $y_1<\cdots< y_m$,  and all joining together at $P$,
where the  gradient $u_x$ blows up as in (\ref{sf0}).
Moreover, assume that, in addition, 
$n$ shocks merge together exactly at $P$.

If $0<\bar t<T$, we then define the {\bf singularity index} of the point $P$ to be the integer
\bel{NP1} N(P)~\doteq~\left\{ \bega{rl} n-2\quad &\hbox{if} ~~m=0,\\[1mm]
m+n-1\quad &\hbox{if} ~~m>0.\enda\right.\eeq

If $\bar t=T$, the {\bf singularity index} of the point $P$ is defined
 as
 \bel{NP2} N(P)~\doteq~\left\{ \bega{rl} n-1\quad &\hbox{if} ~~m=0,\\[1mm]
m+n\quad &\hbox{if} ~~m>0.\enda\right.\eeq
\v
For a solution of (\ref{cblaw}) with initial data $\bar u\in \M$, we set
\bel{NGU}
N(\bar u)~\doteq~\sum_{P_i\in\S} N(P_i).\eeq
\end{definition}

\begin{figure}[ht]
\centerline{\hbox{\includegraphics[width=10cm]{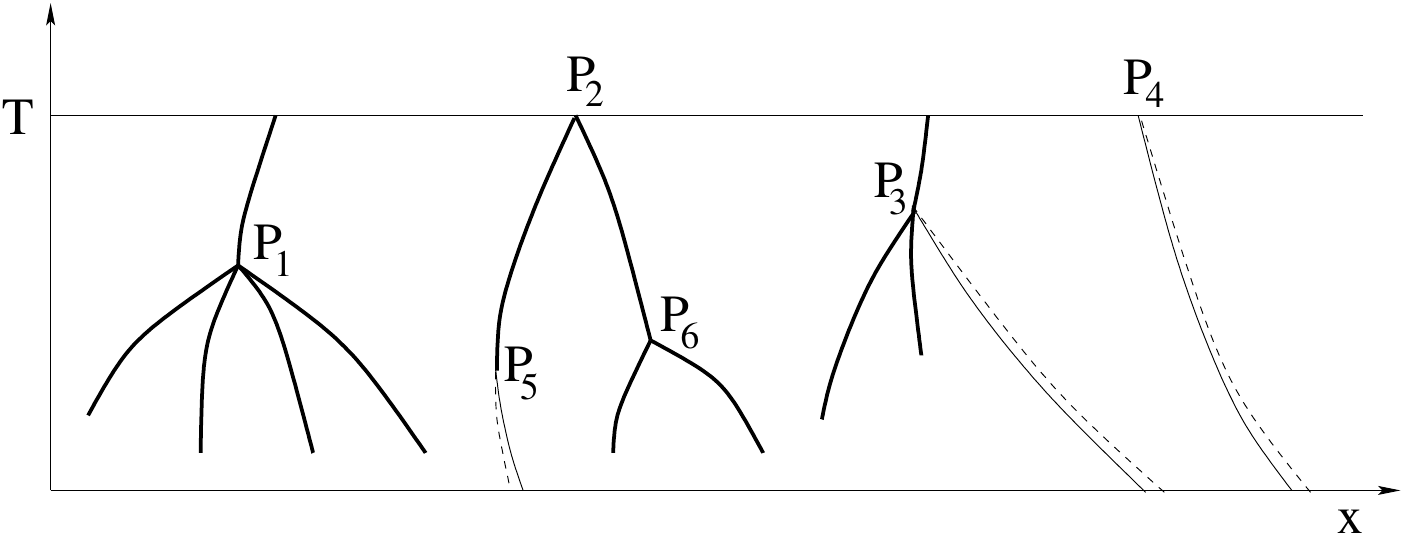}}}
\caption{\small The  singularities at the points $P_1,\ldots,P4$ are non generic. According to Definition~\ref{def:31},
their indices are: $N(P_1) = 2$, $ N(P_2)=1$, $N(P_3)=2$, $N(P_4)=1$.
On the other hand, the singularities at $P_5,P_6$ are generic, hence
$N(P_5)= N(P_6)=0$.  Here the dotted lines indicate characteristics where 
the gradient $u_x$ blows up at the terminal point.
 }
\label{f:gen83}
\end{figure}
%
%

The key step toward a proof of Theorem~\ref{t:2} is the following lemma.

\begin{lemma}\label{l:31}
Let $\bar u\in \M$.  For a fixed $\nu\geq 1$, consider the set 
$\S$  of singular points defined at  (\ref{SP}).  
If $N(\bar u)>0$, then for any $\ve>0$ there exists
$\bar v\in \M$ such that 
\bel{31}
\|\bar v-\bar u\|_{\C^3}\,<\,\ve,\qquad\qquad N(\bar v)\,\leq \,N(\bar u)-1.\eeq
\end{lemma}

Assuming that the lemma is true, we show how it implies the main theorem.

{\bf Proof of Theorem~\ref{t:2}.}
Let $\bar u\in \M$, $\nu\geq 1$  and $\ve_0>0$ be given.  If $N(\bar u) =0$, 
the corresponding solution $u=u(t,x)$ already satisfies all four 
generic conditions (i)--(iv) in Theorem~\ref{t:2}, restricted to the domain
$(t,x)\in [0,T]\times [-\nu,\nu]$.

Next, assume $N(\bar u)>0$. Applying Lemma~\ref{l:31} with  $\ve\doteq \ve_0  /N(\bar u)$, we obtain an initial data $\bar u_1\in \M$ with 
$$\|\bar u_1-\bar u\|_{\C^3}~\leq~\ve,\qquad\qquad N(\bar u_1)\leq N(\bar u)-1.$$
If $N(\bar u_1)=0$, we stop.   Otherwise we continue by induction, constructing a sequence
of initial data $\bar u_1, \bar u_2,\bar u_3,\ldots$ such that
$$\|\bar u_j-\bar u\|_{\C^3}~\leq~j \ve ,\qquad\qquad N(\bar u_j)\leq N(\bar u)-j$$
for all $j\geq 1$.
After a number $r\leq N(\bar u)$ of steps, we obtain an initial data
$\bar u_r$ such that 
$$\|\bar u_r-\bar u\|_{\C^3}~\leq~r \ve \leq\ve_0,\qquad\qquad N(\bar u_r)~=~0.$$

Next, consider the 
set $\M_\nu'\subset\C^3(\R)$ of all initial data $\bar u$ such that the solution 
to (\ref{cblaw}) has finitely many shocks and satisfies all four conditions (i)--(iv), restricted to the domain where $(t,x)\in [0,T]\times [-\nu,\nu]$.

Since $\ve_0>0$ was arbitrary, the previous argument shows that the set $\M_\nu'$ is dense in $\C^3(\R)$.  

We now observe that the conditions (i)--(iv) imply the structural stability
of the solution. Namely, if $\bar u\in \M_\nu'$, then, for any initial data $\bar v\in \C^3(\R)$ 
sufficiently close to $\bar u$, the corresponding solution $v=v(t,x)$ 
will have the same number of shocks as $u$, interacting two at a time, within a neighborhood of $[0,T]\times [-\nu,\nu]$.  We thus conclude that the  set $\M'_\nu$ is also open.

Defining 
$$\M'~\doteq~\M\cap\left(\bigcap_{\nu\geq 1} \M'_\nu\right),$$
we obtain the desired $\G_\delta$ set for which the conclusion of Theorem~\ref{t:2} is satisfied.
\endproof

\v
{\bf Proof of Lemma~\ref{l:31}.}

{\bf 1.} Let $\ve>0$ and an initial data $\bar{u}\in \M$ be given.  
If $N(\bar u)>0$, let $\tau\in \,]0,T]$ be the largest time at which a non-generic interaction 
takes place.   Say, at a point $P=(\tau, \bar x)\in [0,T]\times [-\nu,\nu]$, with $N(P)>0$.

We first consider the case where at least one 
shock curve, say  $t\mapsto x_1(t)$, reaches $P$ at time $\tau$.

\begin{figure}[ht]
\centerline{\hbox{\includegraphics[width=8cm]{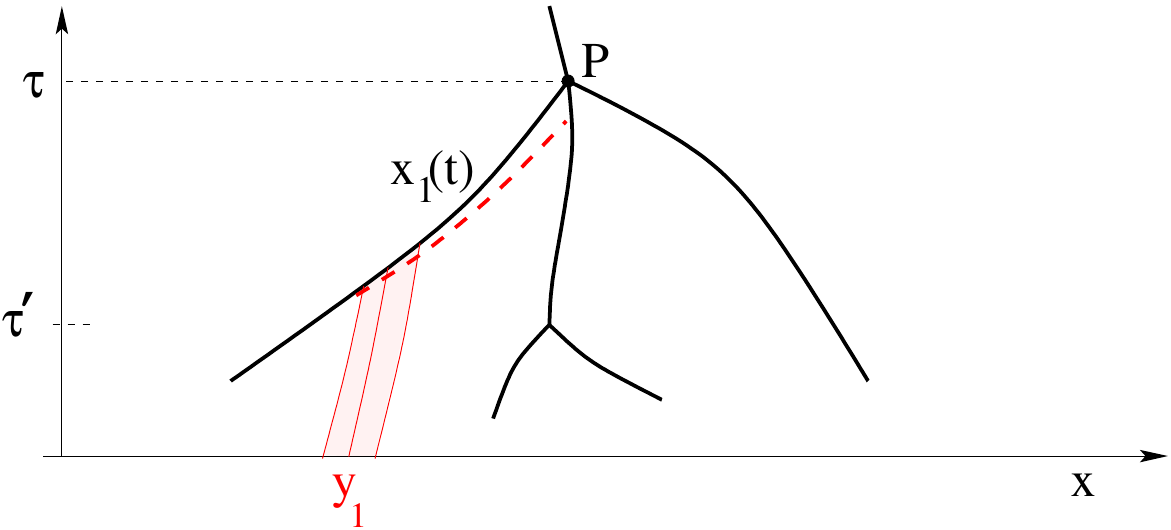}}}
\caption{\small Perturbing the initial data to avoid a triple shock interaction at $P$.
Here the initial data $\bar u$ is slightly increased in a neighborhood of $y_1$, as in 
(\ref{pert3}).
This perturbation produces a forward shift of the shock at $x_1(t)$.
As a result, the three incoming shocks no longer interact all at the same point. 
 }
\label{f:gen79}
\end{figure}

Let $]\tau', \tau[\,$ 
be an open time  interval where no interaction occurs, and
consider a characteristic curve $t\mapsto \xi(t,y_1)$ impinging on the shock 
$x_1(\cdot)$ at some time  $t_1\in \,]\tau', \tau[\,$.  By choosing $\delta>0$ small 
enough, all characteristics starting within the interval $[y_1-\delta, y_1+\delta]$
still impinge on the shock $x_1(\cdot)$ at some time $t \in\,]\tau', \tau[\,$.

We then define the family of perturbed initial data
\bel{pert3}\bar{u}^{\lambda}(x)~\doteq~\bar{u}(x) +\lambda\cdot \vp\left(x-y_1\over \delta\right),
\eeq
where $\vp\in \C^\infty_c(\R)$ is a smooth function such that 
$$
\left\{ \bega{rl} \vp(x) >0\quad\hbox{if} \quad |x|<1,\\[2mm]
\vp(x)=0\quad\hbox{if}\quad |x|\geq 1.\enda\right.
$$

{\bf 2.}
Let $u^{\lambda}$ be  the solution of  the balance law (\ref{cblaw}) with initial data $\bar{u}^{\lambda}$. Our next goal is to compute the shift in  the
position of the perturbed shock $x^\lambda_{1}(t)$  for $t\in \,]\tau',\tau[\,$.

  According to Theorem~2.2 in \cite{BM2},  on the region where $u$ is smooth, the perturbed solution can be described by
\bel{ulap}
u^\lambda(t,x) ~=~u(t,x)+\lambda\cdot  w(t,x) + o(\lambda), \eeq
where $w$ provides a solution to the linearized equation
\bel{pert1}
\begin{cases}
w_t + f'(u) w_x + f''(u) w u_x ~=~g_u(t,x,u,\alpha) w,\\[4mm]
\ds w(0,x)~=~ \vp\left({x-y_1\over \delta}\right).
\end{cases}
\eeq 
To determine
the shift of the shock, defined as
$$\zeta(t)~\doteq ~\lim_{\lambda\to 0+} {x^\lambda_{1}(t)-x_1(t)\over \lambda}\,,$$
we consider the left and right limits of the solution $u$ and of its gradient $u_x$, across the shock:
$$u^\pm(t)~\doteq~\lim_{x\to x_1(t)\pm} u(t, x),\qquad\qquad 
u_x^\pm(t)~\doteq~\lim_{x\to x_1(t)\pm} u_x(t, x).$$

As proved in \cite{BM2}, the shift $\zeta$ is then 
obtained by solving  the ODE
\bel{pshift}\bega{rl}
\dot \zeta(t)&\ds=~
{\partial \over \partial u^+} \Lambda\bigl(u^+(t), u^-(t)\bigr)\cdot \Bigl[w\bigl(t,x_1(t)+\bigr) + \zeta(t) u_x^+(t)\Big] \\[4mm]
&\qquad\qquad\ds
+{\partial \over \partial u^-}\Lambda\bigl(u^+(t), u^-(t)\bigr)\cdot \Bigl[w\bigl(t,x_1(t)-\bigr) + \zeta(t) u_x^-(t)\Big],\enda
\eeq
with initial data $\zeta(\tau')=0$. 
Here $\Lambda$ denotes the Rankine-Hugoniot speed of the shock:
$$
\Lambda(u^+,u^-)~\doteq~{f(u^+) - f(u^-)\over u^+-u^-}~=~\int_{0}^{1}f'(u^-+s(u^+-u^-))ds.
$$
Since $w(0,x)= {\vp\left({x-y_1\over \delta}\right)\geq 0}$, one derives from  (\ref{pert1}) that 
\bel{non-w}
w(t,x)~\geq~0,\qquad (t,x)\in [0,T]\times \R.
\eeq
Moreover, for every $y\in {\,]y_1-\delta,\,y_1+\delta[\,}$, let $t_y$ be the time when the characteristic $t\mapsto \xi(t,y)$ hits the shock $x_1(\cdot)$. For a.e. $t\in [0,t_y]$
one has 
\[
{d\over dt}w(t,\xi(t,y))~=~G_0(t,y)\cdot w(t,\xi(t,y)),
\]
with 
\[
G_0(t,y)~\doteq~g_u\Big(t,\xi(t,y),u\bigl(t,\xi(t,y)\bigr),\alpha(t)\Big) -f''\Big(u\bigl(t,\xi(t,y)
\bigr)\Big) 
u_x\bigl(t,\xi(t,y)
\bigr).
\] 
Hence,
\[
w\bigl(t,\xi(t,y)\bigr)~=~\exp\left({\int_{0}^{t}G_0(s,y)ds}\right)\cdot \vp\left({x-y_1\over \delta}\right)~>~0,\qquad t\in [0,t_y],
\]
and therefore
\bel{sign}
w\bigl(t,x_1(t)+\bigr)~>~0,\qquad {t\in\, \bigl]t_{y_1-\delta},t_{y_1+\delta}\bigr[}\,.
\eeq
By the strictly convexity of $f$ we have
\bel{pala}\left\{
\bega{rl}
\ds\partial_{u^+}\Lambda(u^+,u^-)&\ds=~\int_{0}^{1}s\cdot f''\bigl(u^-+s(u^+-u^-)\bigr)ds~>~0,\\[2mm]
\ds \partial_{u^+}\Lambda(u^+,u^-)&\ds=~\int_{0}^{1}(1-s)\cdot f''\bigl(u^-+s(u^+-u^-)\bigr)ds~>~0,
\enda\right.
\eeq
and (\ref{pshift})-(\ref{non-w}) yield
\bel{ODE4}
\dot \zeta(t)~\geq~G_1(t)\cdot \zeta(t)+{\partial \over \partial u^+} 
\Lambda\bigl(u^+(t), u^-(t)\bigr)\cdot w\bigl(t,x_1(t)+\bigr) ,
\eeq
with 
\[
G_1(t)~=~{\partial \over \partial u^+}\Lambda\bigl(u^+(t), u^-(t)\bigr)\cdot u_x^+(t)
+{\partial \over \partial u^-} 
\Lambda\bigl(u^+(t), u^-(t)\bigr)\cdot u_x^-(t).
\]
Integrating the ODE (\ref{ODE4}) with initial value $\zeta(\tau')=0$, in view of 
(\ref{sign})-(\ref{pala}) we obtain
\[
\zeta(\tau)\,>\,0.
\]
The above computation shows that, by perturbing the initial data in a neighborhood of
the point $y_1$, we can shift the position of the shock $x_1(\tau)$ forward
(taking $\lambda>0$) or backward (taking $\lambda<0$).
By doing so, we can arrange so that the shock $x_1$ hits one of the other shocks
(or one of the characteristics reaching $P$ along which the gradient  $u_x$ blows up)
at a time slightly earlier than $\tau$.   This will decrease the index $N(P)$ at least by one, as required.

By choosing $\bar v = \bar u^\lambda$, with $|\lambda|$ small enough, 
both requirements in (\ref{31}) are satisfied.  
\v
{\bf 3.}
It remains to consider the case where $P= \bigl(T,\xi(T,\bar y)\bigr)$ is a point 
along a characteristic where a new shock is formed exactly at time $T$ 
(see for example the point $P_4$ in Fig.~\ref{f:gen83}). In this case, there exists a neighborhood $\mathcal{V}$ of $\bar{y}$ such that 
\bel{cond2}
\xi_y(t,y)~\geq~0\qquad\forall (t,y)\in [0,T]\times \mathcal{V}.
\eeq
Moreover, a strict inequality holds at all points $(t,y)\not=(T, \bar y)$.

Consider a perturbation of $\bar{u}$ of the form
\[
\bar{u}^{\lambda}(y)~=~\bar{u}(\bar{y})+\lambda\cdot \eta\left({y-\bar{y}\over \delta}\right)\cdot (y-\bar{y}).
\]
Here $\eta$ is the cutoff function at (\ref{cutoff}), and $\delta>0$
is chosen small enough so that the perturbation does not affect any other shock in the solution $u$.
 
Let $\xi^\lambda, v^\lambda$ be the corresponding solutions of (\ref{ODE1}), with $\bar u$ replaced by $\bar u^\lambda$. According to Step~5 in the proof of Theorem \ref{t:1}, the maps  
$$t\mapsto \big(\xi_y(t,\bar{y}), v_y(t,\bar{y})\big),\qquad\quad 
t\mapsto \bigl(\partial_{\lambda}\xi_y^{\lambda}(t,\bar{y}),~\partial_{\lambda}v^{\lambda}_y(t,\bar{y})\bigr),$$
provide a solution to  the system of ODEs
\[
\begin{cases}
\dot X(t)~=~a(t)\cdot Y(t),\\[2mm]
\dot Y(t)~=~b(t)\cdot X(t)+c(t)\cdot Y(t),
\end{cases}
\]
with initial data 
\[\left\{ \bega{l} 
X(0)~=~\xi_y(0,\bar{y})~=~1,\\[1mm]
Y(0)~=~v_y(0,\bar{y})~=~\bar u_x(\bar{y})+\lambda,\enda\right.
\qquad\hbox{and}\qquad 
\left\{\bega{l}
X(0)~=~\partial_{\lambda}\xi_y^{\lambda}(0,\bar{y})~=~0,\\[1mm]
Y(0)~=~\lambda_yv_y^{\lambda}(0,\bar{y})~=~1,\enda\right.
\]
respectively.
Since $\xi_y(T,\bar{y})=0$, we have 
\[
\partial_{\lambda}\xi_y^{\lambda}(T,\bar{y})\bigg|_{\lambda=0}~=~\beta~\neq~0.
\]
Without loss of generality, assume that $\beta>0$. By  continuity, there exist $\delta_1,\delta_2>0$ such that 
%
such that 
\bel{cond3}
\partial_{\lambda}\xi_y^{\lambda}(t,y)~\geq~{\beta\over 2}~>~0\qquad \qquad \forall ~|\lambda|<\delta_1,~~(t,y)\in \Gamma,\eeq
where $\Gamma$ is the rectangle
$$\Gamma~\doteq~\bigl\{ (t,y)\,;~~|t-T|\leq \delta_2, ~~|y-\bar y|\leq \delta_2\bigr\}.$$
On the other hand, by (\ref{cond2}) there exists $\delta_3\in \,]0, \delta_1]$ such that 
\bel{cond4}
\xi^\lambda_y(t,y)\,>\,0\qquad\qquad\hbox{whenever}~~ |\lambda|<\delta_3, ~~~(t,y)\in \Big([0,T]\times \mathcal{V}\Big)\setminus \Gamma.
\eeq
Combining (\ref{cond3}) with (\ref{cond4}) we conclude that, for $0<\lambda<\delta_3$,
the solution to (\ref{ODE1}) satisfies
$$\xi_y^\lambda(t,y)~>~0\qquad\qquad\forall (t,y)\in [0,T]\times \mathcal{V}.$$
By setting $\bar v = \bar u^\lambda$, we thus achieve the second inequality in (\ref{31}).
The first inequality is easily achieved by further shrinking the value of $\lambda>0$.
This completes the proof of Lemma~\ref{l:31}. 
\endproof

\section{Generic behavior of optimally controlled balance laws}
\setcounter{equation}{0}
\label{s:4}

In this section we give a proof of Proposition~\ref{p:11}, constructing an example
where the optimal solution to a controlled balance law always has two shocks interacting exactly at the terminal time $T$.

\begin{figure}[ht]
\centerline{\hbox{\includegraphics[width=11cm]{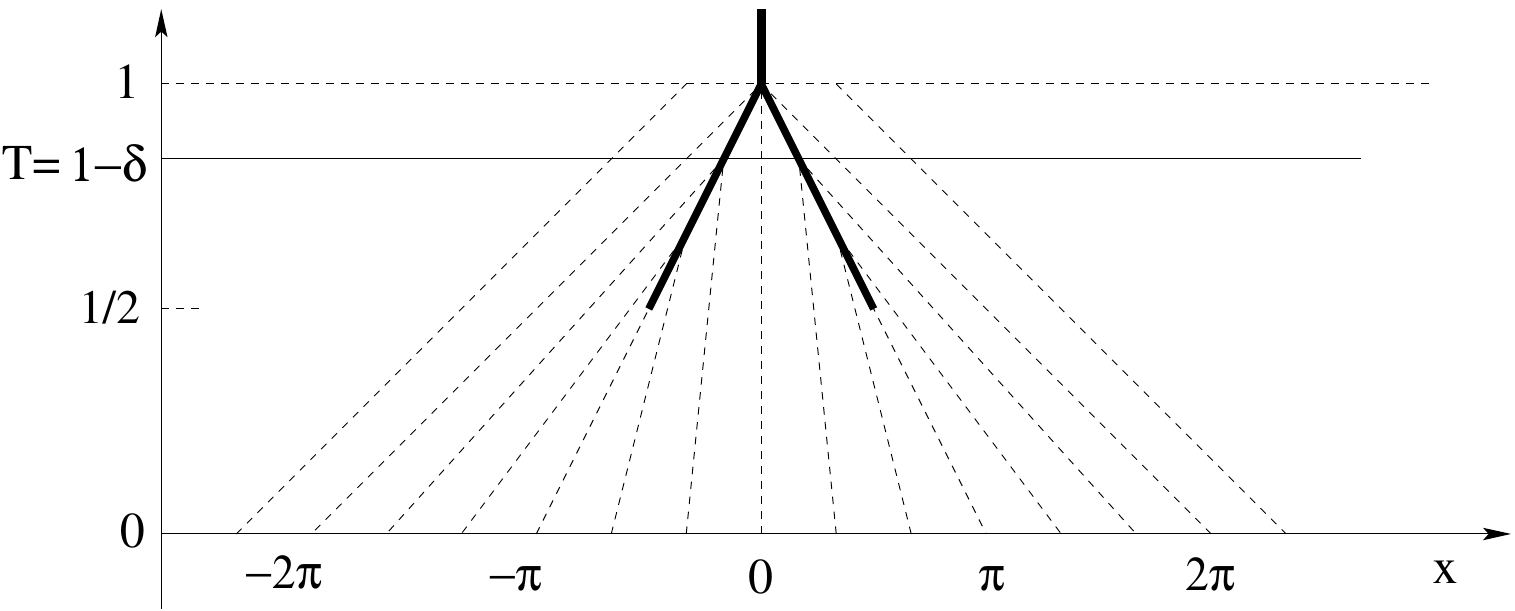}}}
\caption{\small  In the solution to Burgers' equation with initial data (\ref{Bid}), 
two shocks form at time $t=1/2$ and interact at time $t=1$.}
\label{f:gen80}
\end{figure}

As a preliminary, consider the solution of Burgers' equation
\bel{Bur}u_t + \left(u^2\over 2\right)_x~=~0,\eeq
with initial data
\bel{Bid}
u(0,x)~=~\bar u(x)~\doteq~\left\{\bega{cl} \sin x-x\quad &\hbox{if}\quad x\in [-2\pi,2\pi],\\[1mm]
2\pi\quad &\hbox{if}\quad x\leq -2\pi,\\[1mm]
-2\pi\quad &\hbox{if}\quad x>2\pi.\enda\right.
\eeq
Notice that the gradient $\bar u_x$
has two negative  minima at the points
$$x_1=-\pi, \quad x_2=\pi,\qquad\hbox{with} \quad {\bar u_x(x_1)\,=\,\bar u_(x_2)\,=\,-2.}$$
In the solution of (\ref{Bur})-(\ref{Bid}), new shocks are formed at time $t=1/2$, 
at the points
$p_1=-\pi/2$, $p_2=\pi/2$.
As shown in Fig.~\ref{f:gen80},
these shocks are located at the points
$$x_1(t)~=~\pi(t-1) , \qquad\qquad  x_2(t) ~=~ \pi(1-t),\qquad\qquad t\in [1/2, ~1]\,,$$
and merge together at time $t=1$, at $x=0$.

We now consider an optimal control problem of the form  (\ref{min1}), on the time interval 
$t\in [0,T]$,  with $T=1-\delta$ for some $\delta>0$ small, with dynamics
\bel{Burs}u_t + \left(u^2\over 2\right)_x~=~-\eta(x)\,\alpha(t).\eeq
and initial data (\ref{Bid}).
Here $\eta$ is a smooth odd function with compact support such that 
\bel{edef}
\eta(x)~=~\left\{ \bega{rl} x\quad \hbox{if} ~~|x|\leq 2\pi,\\[1mm]
0\quad \hbox{if} ~~|x|\geq 4\pi,\enda\right.\qquad\qquad \eta(-x) \,=\, -\eta(x).\eeq
We observe that, since both the initial condition $\bar u$ in (\ref{Bid}) and the source term in (\ref{Burs}) are odd functions of $x$, the solution $u=u(t,x)$ will also 
satisfy $u(t,-x) =-u(t,x)$ for all $t,x$.  This symmetry simplifies some of the 
computations.

As running cost and terminal cost functions, we choose
\bel{cost2}\Phi(\alpha)\,=\,\alpha^2,\qquad\qquad 
\Psi(u)\,=\,\left\{ \bega{cl}  (u-1)^2(u+1)^2\qquad &\hbox{if}\quad  |u|\leq 1,\\[2mm]
0\qquad &\hbox{if} \quad |u|>1.\enda\right.\eeq

\begin{figure}[ht]
\centerline{\hbox{\includegraphics[width=11cm]{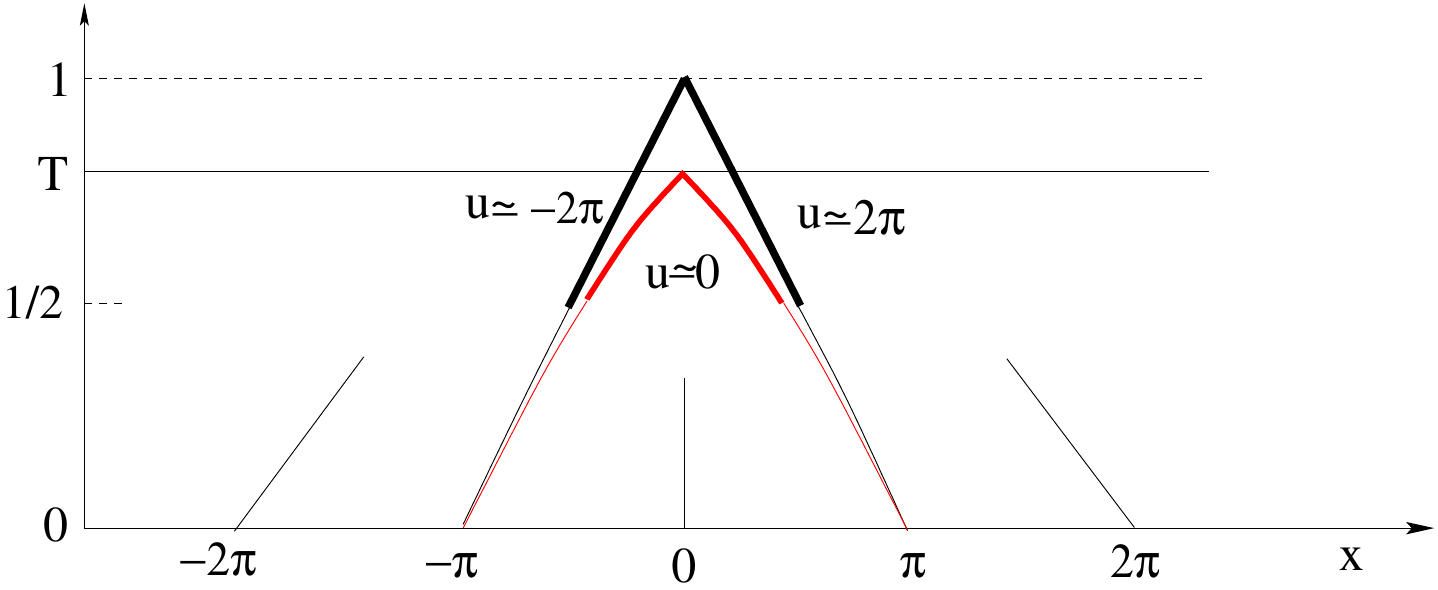}}}
\caption{\small When $t\approx 1$, the solution $u(t,\cdot)$ takes values
close to zero between the two shocks. If no control is applied, the two shocks 
interact at time $t=1$.  To reduce the terminal cost,
it is convenient to implement a control which makes the two shocks merge together within time $T$.}
\label{f:gen81}
\end{figure}
The motivation behind these choices is explained in Fig.~\ref{f:gen81}.
The terminal cost $\Psi$ penalizes values contained in  $[-1,1]$.  In particular,
calling $x_1(t)<x_2(t)$ the positions of the two shocks,
the terminal cost function will be of the form
\bel{IT}\int_{x_1(T)}^{x_2(T)} \Psi\bigl(u(T,x)\bigr)\, dx~\approx~\int_{x_1(T)}^{x_2(T)} \Psi(0)\, dx~=~
x_2(T)-x_1(T).\eeq
This makes it more convenient to have the two shocks merge together within 
time $T$.
Choosing a control $\alpha(t)>0$ on the right hand side of (\ref{Burs}), the values of $u(t,x)$
increase for $x<0$ and decrease for $x>0$.  This increases the speed of the shock at $x_1(\cdot)$
and decreases the speed of the shock at $x_2(\cdot)$.    As a result, the two shocks 
will interact
at an earlier time.   Since there is a running cost penalizing large control values,
the size of the control should be chosen as small as possible: just enough to
make the two shocks interact exactly at the terminal time $T$. There is no advantage in
making the shocks interact at any earlier time.

The proof of Proposition~\ref{p:11} will be given in three steps.
\v
{\bf 1.} Choosing $\delta>0$ small 
enough, 
we claim that every optimal solution will contain exactly two shocks.

Indeed, consider the trivial control $\alpha(t)\equiv 0$.
By (\ref{IT}),  when $T=1-\delta$ the terminal cost will have size 
$$\O(1)\cdot \bigl(x_2(T)\bigr)-x_1(T)\bigr)~=~\O(1)\cdot \delta.$$
Letting $\delta\to 0+$ the total cost associated with the zero control
thus approaches zero.  
If now $\alpha^*(\cdot)$ is an optimal control, its total cost cannot be any larger.
In particular, the running cost of $\alpha^*(\cdot)$  
must also approach zero.  Therefore
\bel{b-control}
\|\alpha^*\|_{{\bf L}^1}~=~\O(1)\cdot \|\alpha^*\|_{{\bf L}^2}~=~\O(1)\cdot\left( \int_0^T
\bigl|\alpha^*(t)\bigr|^2\, dt\right)^{1/2} ~=~\O(1)\cdot \sqrt \delta.
\eeq
 It now remains to show that, for the balance law (\ref{Burs}),  a control function $\alpha(\cdot)$ with sufficiently small 
 $\L^1$ norm
cannot produce additional shocks within the interval $t\in [0,1]$.

Toward this goal, given a control $\alpha(\cdot)$, let $\bigl(\xi^{\alpha}(t,y),v^{\alpha}(t,y)\bigr)$ be the solution to the system of ODEs for characteristics:
\bel{ODE1-pB}
\begin{cases}
\dot{\xi}~=~v,\\
\dot{v}~=~-\eta(\xi)\alpha,
\end{cases}
\qquad\mathrm{with}\qquad\quad
\begin{cases}
\xi(0)~=~y,\\
v(0)~=~\bar{u}(y).
\end{cases}
\eeq
Calling $\xi=\xi(t,y)$ the solution to (\ref{ODE1-pB})  when the control is identically zero,
one checks that
\bel{per-xi}
\|\xi^{\alpha}-\xi\|_{\mathcal{C}^3([0,T]\times \R)}~\leq~\O(1)\cdot \|\alpha\|_{\L^1}
~=~\O(1)\cdot \sqrt\delta\,.
\eeq
In the case of zero control, we compute
\bel{xi-y}
 \xi_y(t,y)~=~1+t\,\bar u_x(y)~=~\left\{ \bega{cl}  1+ t\, (\cos y-1)\quad  &\hbox{if} \quad |y|\leq  2\pi,\\[2mm]
 1\quad &\hbox{if}  \quad|y|> 2\pi.
\enda\right.
\eeq
Hence $\xi_y(t,y)$ can only vanish for $y\in \,]-2\pi,0[\,\cup \,]0,2\pi[\,$. In this case, the time $T(y)$ at which $\xi_y(t,y)$ vanishes is given by 
\bel{Ty}
T(y)~=~{1\over 1-\cos (y)},\qquad y\in \,]-2\pi,0[\,\cup \,]0,2\pi[\,.\eeq
By (\ref{xi-y}) and (\ref{per-xi}) it follows that, for $\delta>0$ sufficiently small,  the gradient
$\xi^{\alpha}_y(t,y)$ can vanish at some time $t\in [0,1]$ only if $|y|\in [\pi/4,7\pi/4]$. Moreover, by (\ref{per-xi}) we have
\[
\bigl|\partial_{t}\xi^{\alpha}_y(t,y)\bigr|\,\geq\,\bigl|\partial_{t}\xi_y(t,y)\bigr|-\O(1)\cdot \sqrt \delta\,=\,|1-\cos y|-\O(1)\cdot \sqrt\delta\,\geq\,{1\over 4}\qquad
\hbox{for}~~ |y|\in \left[{\pi\over 4}\,,\,{7\pi\over 4}\right].
\]
 Let $T^{\alpha}(y)$ be such that $ \xi^{\alpha}_{y}(T^{\alpha}(y),y)=0$. By the implicit function theorem we obtain 
 \[
 \|T^{\alpha}-T\|_{\mathcal{C}^2\bigl([-7\pi/4,-\pi/4]\cup[\pi/4,7\pi/4]\bigr)}~=~\O(1)\cdot \sqrt\delta.
 \]
 Differentiating (\ref{Ty}) one obtains
 \[
 {d^2\over dy^2} T(y)~=~{1+\sin^2y-\cos y\over (1-\cos y)^3}~\geq~1-{\sqrt{2}\over 2}\qquad \hbox{whenever}\quad |y|\in \left[{\pi\over 4}\,,\,{7\pi\over 4}\right].
 \]
 As a consequence, if $\|\alpha\|_{\L^1} =\O(1)\cdot \sqrt\delta>0$ is sufficiently small, then
 the function $T^{\alpha}(\cdot)$  attains a local minimum at  only two points: $y^{\delta}_1$ and $y^{\delta}_2$, with 
 \bel{y12d} y_1^\delta\,\to \,-\pi,\qquad y_2^\delta \,\to\, \pi\qquad\hbox{as}\quad \delta\to 0.\eeq
%
%
In the remainder of the proof, we can thus assume that the optimal solution is piecewise
continuous, with exactly two shocks.   Assume that  these shocks interact at a  time $\tau^*\not=T$. In the next two steps we prove the followings: 
\begin{itemize}
\item If $\tau^*>T$ then  to reduce the terminal cost, it is convenient to increase the size of the control,
so that the distance $x_2(T)-x_1(T)$ between the two shocks at the terminal time 
is decreased.
\item  Otherwise, if $\tau^*<T$ then it is convenient to reduce the size of the control,
letting the two shocks interact at a later time.
\end{itemize}
%

{\bf 2.} Let the pair $(\alpha^*,u^*)$ denote an  optimal control and an optimal solution. Assume that the two shocks in $u^*$, located at  $x_1^*(t)<0<x^*_2(t)$, which are formed at 
time $t^*$, remain separated for all $t\in [t^*,T]$. In particular,
$x_1^*(T) < x_2^*(T)$.

We now choose a time $\tau\in \,]t^*, T[\,$ after the shocks have formed, and 
consider a perturbation of the control defined by
\bel{aep}\alpha^\ve(t)~\doteq ~\left\{ \bega{cl} \alpha^*(t)\quad&\hbox{if} ~~t\in [0, \tau],
\\[1mm] 
 \alpha^*(t)+\ve\quad&\hbox{if} ~~t\in \,]\tau,T].\enda\right.\eeq
We claim that, for $\ve>0$ small, the new control yields a lower cost.

Indeed, let 
$$u^\ve~=~ u+\ve w + o(\ve)$$ be the solution to (\ref{Burs})
corresponding to the perturbed control.  By (\ref{aep}) we trivially have $u^\ve = u$ 
for $t\in [0, \tau]$. For $t\in [\tau,T]$,
outside the shocks the leading order perturbation $w$ is determined by the linearized system
\bel{pert2}
w_t + u^* w_x +   u^*_x w ~=~-\eta(x),\qquad\qquad  w(\tau,x)~=~ 0.
\eeq 
Similarly, call
$$x_1^\ve(t) ~=~x_1^*(t) + \ve {\zeta_1(t)} + o(\ve),\qquad  x_2^\ve(t) ~=~x_2^*(t)
 + \ve {\zeta_2(t)} + o(\ve)$$
the positions of the two shocks in the perturbed solution.   By symmetry,
$\zeta_1(t)= - \zeta_2(t)$.  The shifts $\zeta_1,\zeta_2$ are 
then determined by solving the ODEs
\bel{pshift}
\dot \zeta_i(t)~=~ {1\over 2} \Big[w\bigl(t,x^*_i(t)-\bigr)+w\bigl(t,x^*_i(t)+\bigr)\Big]
 + {1\over 2}  \Big[u^*_x\bigl(t,x^*_i(t)-\bigr)+u^*_x\bigl(t,x^*_i(t)+\bigr)\Big] \cdot \zeta_i(t), \eeq
 with initial data
 \bel{ishift}\zeta_1(\tau)\, =\, \zeta_2(\tau)\,=\,0.\eeq
 
 Consider an initial point $y\in [0, 2\pi]$ and
 let $t\mapsto (\xi^*(t,y)),v^*(t,y))$ be the solution of the Cauchy problem (\ref{ODE1-pB}) with $\alpha=\alpha^*$. Before the characteristic $t\mapsto \xi^*(t,y)$ reaches the shock at $x_2^*(t)$,  the perturbation $w$   satisfies
 \[
 {d\over dt}w(t,\xi^*(t,y))+u^*_x(t,\xi^*(t,y)) \cdot w(t,\xi^*(t,y))~=\,-\eta\bigl(\xi^*(t,y)\bigr),\qquad\quad w\bigl(\tau,\xi^*(t,y)\bigr)~=~0.
 \]
Hence, also in view of (\ref{edef}), there exists a constant $c_0>0$ such that 
 \bel{per1}
 w\bigl(t,\xi^*(t,y)\bigr)~\leq~-c_0\cdot \int_\tau^{t}\xi^*(s,y)ds\,. \eeq
In turn, this implies 
\[
w(t,x^*_2(t)-)+w(t,x^*_2(t)+)~\leq~-c_1,
\]
 for some constant $c_1>0$ and all $t\in \bigl[(\tau+T)/2, ~T\bigr]$. Thus,  from (\ref{pshift}), we obtain
 \[
 \zeta_2(T)\,=\,- c_2,\qquad\qquad \zeta_1(T) \,=\, c_2\,,
 \]
 for some constant $c_2>0$.
For $\ve>0$ sufficiently small,  the cost of the perturbed control satisfies
\[
\bega{l}
J[\alpha^\ve]-J[\alpha^*]\\[4mm]
\qquad \ds= ~\int_{\tau}^T 2\ve \,\alpha^*(t)\, dt  -
\left( \int_{x_1^*(T)}^{x_1^\ve(T)} + \int^{x_2^*(T)}_{x_2^\ve(T)}\right)
 \bigl(u^*(T,x) -1\bigr)^2\bigl(u^*(T,x)+1\bigr)^2\, dx + \O(1)\cdot\ve^2
\\[4mm]
\qquad \leq ~{2\ve} \cdot\|\alpha^*\|_{\L^1}  + \ve\, \bigl( {\zeta_2(T)} - {\zeta_1(T)}\bigr)+ \O(1)\cdot \ve^2 + \O(1)\cdot  \ve \bigl(x_2^*(T)- x_1^*(T)\bigr)\\[4mm]
\qquad \leq ~\ve \cdot\Big( \|\alpha^*\|_{\L^1}  -2c_2 +\O(1)\cdot \sqrt\delta
\Big) + \O(1)\cdot \ve^2 ~<~0,\enda
\]
provided that $\ve,\delta>0$ are small enough.
This yields a contradiction. 

\v
{\bf 3.} Next, consider the case where in the optimal solution the two shocks merge at a time $\tau^*<T$. 
As a preliminary we observe that, when the control $\alpha(t)\equiv 0$ is identically zero, all characteristics $t\mapsto \xi(t,y)$ starting at a point $y\in [-\pi,\pi]$ impinge on one of the two shocks (see Fig.~\ref{f:gen80}). 
Along all other characteristics one has
$$\bigl|u(t,\xi(t,y))\bigr|~=~\bigl|\bar u(y) \bigr| ~\geq~\pi.$$
Next, let $\alpha(\cdot)$ be a control, with $\|\alpha\|_{\L^1}=\O(1)\cdot\sqrt\delta$ 
sufficiently small, such that the corresponding solution $u$ of (\ref{Burs}), (\ref{Bid})
has two shocks interacting at time $\tau^*<T$. In this case, the only characteristics that 
do not impinge on one of the shocks within time $\tau^*$ must originate at a point $y\leq -\pi$ or 
$y\geq \pi$. By continuity, we have
$$\Big|u\bigl( \xi(T,y)\bigr)\Big| ~=~\Big|\bar u(y) + \O(1)\cdot \sqrt\delta\Big|~\geq~\pi - \O(1)\cdot \sqrt\delta.$$
In view of  (\ref{cost2}) this yields
$$\Psi\bigl(u( \xi(T,y))\bigr) ~=~0.$$
Hence for this solution the terminal cost is zero, and only the running cost related to the control $\alpha(\cdot)$ is present.   In particular, we have
\[
J[\alpha^*]~=~\int_{0}^{T}[\alpha^{*}(t)]^2dt.
\]
We now consider the perturbed controls
$$\alpha^\ve(t)~\doteq~(1-\ve) \alpha^*(t).$$
By continuous dependence,  for $\ve>0$ small the corresponding solution $u^\ve$ still contains two shocks 
merging before time $T$. Hence the previous arguments yield
\[
J[\alpha^\ve]~=~\int_{0}^{T}[\alpha^\ve(t)]^2\,dt~=~(1-\ve)^2 J[\alpha^*]~<~J[\alpha^*],
\]
reaching a contradiction.
\v
{\bf 4.} In both of the above cases, we concluded that if the shocks interact at a time $\tau\not=T$,
the control $\alpha(\cdot)$ is not optimal. The previous analysis was performed in connection with the balance law (\ref{cblaw})
with 
$$f(u)={u^2\over 2},\qquad \qquad g(t,x,u,\alpha)~\doteq~-\eta(x)\alpha,$$
and $\bar u, \Phi,\Psi$ as in (\ref{Bid}), (\ref{cost2}).
To achieve a proof of Proposition~\ref{p:11}, it now suffices to observe that the previous steps {\bf 1--3} 
remain valid if the dynamics, cost functions and initial data  
are replaced by any perturbations
$(\Tilde f,\Tilde g,\Tilde \Phi,\Tilde \Psi,\Tilde u)$, provided that 
$$\|\Tilde f-f\|_{\C^3}\,,\qquad \|\Tilde g-g\|_{\C^3}\,,\qquad  \|\Tilde \Phi-\Phi\|_{\C^2}\,,
\qquad 
 \|\Tilde \Psi-\Psi\|_{\C^2}\,,\qquad  \|\Tilde u-\bar u\|_{\C^3}\,,$$
are sufficiently small.
\endproof

\v
{\small
{\bf Acknowledgments.} The research of A.\,Bressan was partially supported by NSF with
grant  DMS-2306926, ``Regularity and approximation of solutions to conservation laws".
The research of K.\,T.\,Nguyen was partially supported by NSF with grant DMS-2154201, ``Generic singularities and fine regularity structure for nonlinear partial differential equations".}

\end{document}